
  \magnification 1200
  \input amssym  


  \newcount\fontset
  \fontset=1
  \def \dualfont#1#2#3{\font#1=\ifnum\fontset=1 #2\else#3\fi}

  \dualfont\bbfive{bbm5}{cmbx5}
  \dualfont\bbseven{bbm7}{cmbx7}
  \dualfont\bbten{bbm10}{cmbx10}

  \font \eightbf = cmbx8
  \font \eighti = cmmi8 \skewchar \eighti = '177
  \font \eightit = cmti8
  \font \eightrm = cmr8
  \font \eightsl = cmsl8
  \font \eightsy = cmsy8 \skewchar \eightsy = '60
  \font \eighttt = cmtt8 \hyphenchar\eighttt = -1

  \font \sixi = cmmi6 \skewchar \sixi = '177
  \font \sixrm = cmr6
  \font \sixsy = cmsy6 \skewchar \sixsy = '60
  \font \tensc = cmcsc10
  
  \font \titlefont = cmbx12
  \scriptfont \bffam = \bbseven
  \scriptscriptfont \bffam = \bbfive
  \textfont \bffam = \bbten

  \newskip \ttglue

  \def \eightpoint {\def \rm {\fam0 \eightrm }%
  \textfont0 = \eightrm
  \scriptfont0 = \sixrm \scriptscriptfont0 = \fiverm
  \textfont1 = \eighti
  \scriptfont1 = \sixi \scriptscriptfont1 = \fivei
  \textfont2 = \eightsy
  \scriptfont2 = \sixsy \scriptscriptfont2 = \fivesy
  \textfont3 = \tenex
  \scriptfont3 = \tenex \scriptscriptfont3 = \tenex
  \def \it {\fam \itfam \eightit }%
  \textfont \itfam = \eightit
  \def \sl {\fam \slfam \eightsl }%
  \textfont \slfam = \eightsl
  \def \bf {\fam \bffam \eightbf }%
  \textfont \bffam = \bbseven
  \scriptfont \bffam = \bbfive
  \scriptscriptfont \bffam = \bbfive
  \def \tt {\fam \ttfam \eighttt }%
  \textfont \ttfam = \eighttt
  \tt \ttglue = .5em plus.25em minus.15em
  \normalbaselineskip = 9pt
  \def \MF {{\manual opqr}\-{\manual stuq}}%
  \let \sc = \sixrm
  \let \big = \eightbig
  \setbox \strutbox = \hbox {\vrule height7pt depth2pt width0pt}%
  \normalbaselines \rm }



  \newcount \secno \secno = 0
  \newcount \stno \stno =0
  \newcount \eqcntr \eqcntr=0

  \def \ifn #1{\expandafter \ifx \csname #1\endcsname \relax }

  \def \track #1#2#3{\ifn{#1}\else {\tt\ [#2 \string #3] }\fi}

  \def \advseqnumbering {\global \advance \stno by 1 \global \eqcntr=0}

  \def \current {\number \secno \ifnum \number \stno = 0\else
    .\number \stno \fi }

  \def \laberr#1#2{\message{*** RELABEL CHECKED FALSE for #1 ***}
      RELABEL CHECKED FALSE FOR #1, EXITING.
      \end}

  \def \syslabel#1#2{%
    \ifn {#1}%
      \global \expandafter 
      \edef \csname #1\endcsname {#2}%
    \else
      \edef\aux{\expandafter\csname #1\endcsname}%
      \edef\bux{#2}%
      \ifx \aux \bux \else \laberr{#1=(\aux)=(\bux)} \fi
      \fi
    \track{showlabel}{*}{#1}}

  \def \subeqmark #1 {\global \advance\eqcntr by 1
    \edef\aux{\current.\number\eqcntr}
    \eqno {(\aux)}
    \syslabel{#1}{\aux}}

  \def \eqmark #1 {\advseqnumbering
    \eqno {(\current)}\syslabel{#1}{\current}}

  \def \fcite#1#2{\syslabel{#1}{#2}\lcite{#2}}

  \def \label #1 {\syslabel{#1}{\current}}

  \def \lcite #1{(#1\track{showcit}{$\bullet$}{#1})}

  \def \cite #1{[{\bf #1}\track{showref}{\#}{#1}]}

  \def \scite #1#2{{\rm [\bf #1\track{showref}{\#}{#1}{\rm \hskip 0.7pt:\hskip 2pt #2}\rm]}}


 \def \Headlines #1#2{\nopagenumbers
    \advance \voffset by 2\baselineskip
    \advance \vsize by -\voffset
    \headline {\ifnum \pageno = 1 \hfil
    \else \ifodd \pageno \tensc \hfil \lcase {#1} \hfil \folio
    \else \tensc \folio \hfil \lcase {#2} \hfil
    \fi \fi }}

  \def \Title{\centerline{\titlefont \ucase{\titletextOne}}
    \edef\aux{\titletextTwo}\edef\bux{\null}%
    \ifx\aux\bux \else 
      \smallskip
      \centerline{\titlefont \ucase{\titletextTwo}}
      \fi}

  \def \Date #1 {\footnote {}{\eightit Date: #1.}}


  \def \ucase #1{\edef \auxvar {\uppercase {#1}}\auxvar }
  \def \lcase #1{\edef \auxvar {\lowercase {#1}}\auxvar }

  \def \goodbreak {\vskip0pt plus.1\vsize \penalty -250 \vskip0pt
plus-.1\vsize }

  \def \section #1{\global\def \SectionName{#1}\stno = 0 \global
\advance \secno by 1 \bigskip \bigskip \goodbreak \noindent {\bf
\number \secno .\enspace #1.}\medskip \noindent \ignorespaces}

  \long \def \sysstate #1#2#3{%
    \advseqnumbering
    \medbreak \noindent 
    {\bf \current.\enspace #1.\enspace }{#2#3\vskip 0pt}\medbreak }
  \def \state #1 #2\par {\sysstate {#1}{\sl }{#2}}
  \def \definition #1\par {\sysstate {Definition}{\rm }{#1}}
  \def \remark #1\par {\sysstate {Remark}{\rm }{#1}}


  \def \proof {\medbreak \noindent {\it Proof.\enspace }}
  \def \proofend {\ifmmode \eqno \square \else \hfill \square
\looseness = -1 \medbreak \fi }

  \def \$#1{#1 $$$$ #1}
  \def \=#1{\buildrel \hbox{\sixrm #1} \over =}

  \def \pilar #1{\vrule height #1 width 0pt}

  \def \Item #1{\smallskip \item {{\rm #1}}}
  \newcount \zitemno \zitemno = 0

  \def \izitem {\zitemno = 0}
  \def \zitemplus {\global \advance \zitemno by 1\relax}
  \def \rzitem{\romannumeral \zitemno}
  \def \rzitemplus {\zitemplus \rzitem}
  \def \zitem {\Item {{\rm(\rzitemplus)}}}

  \newcount \nitemno \nitemno = 0
  
  \def \nitem {\global \advance \nitemno by 1 \Item {{\rm(\number\nitemno)}}}

  \newcount \aitemno \aitemno = 0
  \def\boxlet#1{\hbox to 6.5pt{\hfill #1\hfill}}
  
  \def \aitem {\Item {(\ifcase \aitemno \boxlet a\or \boxlet b\or
\boxlet c\or \boxlet d\or \boxlet e\or \boxlet f\or \boxlet g\or
\boxlet h\or \boxlet i\or \boxlet j\or \boxlet k\or \boxlet l\or
\boxlet m\or \boxlet n\or \boxlet o\or \boxlet p\or \boxlet q\or
\boxlet r\or \boxlet s\or \boxlet t\or \boxlet u\or \boxlet v\or
\boxlet w\or \boxlet x\or \boxlet y\or \boxlet z\else zzz\fi)} \global
\advance \aitemno by 1}

  \newcount \footno \footno = 1
  \newcount \halffootno \footno = 1
  \def \footcntr {\global \advance \footno by 1
  \halffootno =\footno
  \divide \halffootno by 2
  $^{\number\halffootno}$}
  \def \fn#1{\footnote{\footcntr}{\eightpoint#1\par}}

  \begingroup
  \catcode `\@=11
  \global\def\eqmatrix#1{\null\,\vcenter{\normalbaselines\m@th%
      \ialign{\hfil$##$\hfil&&\kern 5pt \hfil$##$\hfil\crcr%
	\mathstrut\crcr\noalign{\kern-\baselineskip}%
	#1\crcr\mathstrut\crcr\noalign{\kern-\baselineskip}}}\,}
  \endgroup


  \font\mf=cmex10
  \def\union {\mathop{\raise 9pt \hbox{\mf S}}\limits}
  \def\inters{\mathop{\raise 9pt \hbox{\mf T}}\limits}

  \def \<{\left \langle \vrule width 0pt depth 0pt height 8pt }
  \def \>{\right \rangle }  
  
  \def \ds{\displaystyle}
  \def \and {\hbox {,\quad and \quad }}

  \def \imply {\mathrel{\Rightarrow}}
  \def \for #1{,\quad \forall\,#1}
  \def \square {\hbox {$\sqcap \!\!\!\!\sqcup $}}
  
  \def \stress #1{{\it #1}\/}
  \def \inv {^{-1}}
  \def \*{\otimes}

  \newcount \bibno \bibno =0
  \def \newbib #1{\global\advance\bibno by 1 \edef #1{\number\bibno}}
  \def \bibitem #1#2#3#4{\smallskip \item {[#1]} #2, ``#3'', #4.}
  \def \references {
    \begingroup
    \bigskip \bigskip \goodbreak
    \eightpoint
    \centerline {\tensc References}
    \nobreak \medskip \frenchspacing }

  \font \tensc = cmcsc10
  
  \font\rs=rsfs10

  \def\reg{\lambda}

  \def\c{\subseteq}
  \def\leq{\leqslant}
  \def\geq{\geqslant}
  \def\emptyset{\varnothing}
  \def\bitem{\Item{$\bullet$}}
  \def\qcite#1#2{\cite{#1}}   


  \def\S{{\cal S}}
  
  \def\T{{\cal T}}
  \def\E{E}   
  \def\spec{\widehat \E}
  \def\speci{\spec_\infty}
  \def\specz{\spec}
  \def\tightbox{{\hbox{\eightsl tight}}}
  \def\spect{\spec_\tightbox}
  
  \def\BS{{\cal B}}
  \def\B{{\cal B}}
  \def\I{{\cal I}}
  \def\P{{\cal P}}

  \def\mycalfont{\ifdim\normalbaselineskip = 12pt \rs \else \rssmall \fi}
  
  \def\U{\Omega}

  \def\inf{\wedge}
  \def\sup{\vee}
  \def\nega{\neg\,}  

  \def\alt#1{\hat #1}

  \def\its{\Cap}



  \def \newbib #1#2{\global\advance\bibno by 1 \edef
#1{\number\bibno}}

  \newbib\actions{E}
  \newbib\Kumjian{K}
  \newbib\Lawson{L}
  \newbib\PatBook{P}
  \newbib\Szendrei{S}


  \def\titletextOne{TIGHT REPRESENTATIONS OF SEMILATTICES}
  \def\titletextTwo{\null}
  \def\titletextTwo{AND INVERSE SEMIGROUPS} 

  \Headlines {\titletextOne \ \titletextTwo} {R.~Exel}

  \null\vskip -1cm
  \Title
  \footnote{\null} 
  {\eightrm 2000 \eightsl Mathematics Subject Classification:
  \eightrm 
  20M18,  
  20M30.  
  }

  \bigskip
  \centerline{\tensc 
    R.~Exel\footnote{*}{\eightpoint Partially supported by
CNPq.}}

  \bigskip
  \Date{13 Mar 2007}

  \midinsert 
  \narrower \narrower
  \eightpoint \noindent
  By a \stress{Boolean inverse semigroup} we mean an inverse semigroup
whose semilattice of idempotents is a Boolean algebra.  We study
representations of a given inverse semigroup $\S$ in a Boolean inverse
semigroup which are \stress{tight} in a certain well defined technical
sense.  These representations are supposed to preserve as much as
possible any trace of \stress{Booleannes} present in the semilattice
of idempotents of $\S$.
  After observing that the Vagner--Preston representation is not
tight, we exhibit a canonical tight representation for any inverse
semigroup with zero, called the \stress{regular representation}.  We
then tackle the question as to whether this representation is
faithful, but it turns out that the answer is often negative.  The
lack of faithfulness is however completely understood as long as we
restrict to \stress{continuous} inverse semigroups, a class
generalizing the $E^*$-unitaries.
  \endinsert

\section{Introduction}
  We shall say that an inverse semigroup $\S$ is a \stress{Boolean
inverse semigroup}, if $E(\S)$, the semilattice of idempotents of
$\S$, admits the structure of a Boolean algebra whose order coincides
with the usual order on $E(\S)$.

Boolean inverse semigroups are quite common, a well known  example being the
semigroup $\I(X)$ of all partially defined bijections on $X$.  The
semilattice of idempotents of $\I(X)$ coincides with the Boolean
algebra $\P(X)$ of all subsets of $X$, this being the reason why 
$\I(X)$
is a indeed a Boolean inverse semigroup.

Given an inverse semigroup $\S$ one might like to study how far it is
from being a Boolean inverse semigroup by considering homomorphisms 
  $$
  \sigma:\S\to\B,
  $$
  into some Boolean inverse semigroup $\B$.  Simply requiring $\sigma$
to be a semigroup homomorphism completely sidesteps the issue since,
in case $\S$ itself happens to be a Boolean inverse semigroup, a mere
semigroup homomorphism has no reason to respect the Boolean algebra
structures involved.

To deal with this situation we propose to consider a special class of
homomorphisms called \stress{tight representations} (see
Definition \fcite{DefISGTightRep}{6.1}), which applies to every
inverse semigroup with zero.  In case $\S$ is a Boolean inverse
semigroup we prove in Proposition \fcite{BisgToBisg}{6.2} that tight
representations are precisely those which restrict to a homomorphism
$\sigma: E(\S) \to E(\B)$ in the category of Boolean algebras.

One of the most important homomorphisms from an inverse semigroup $\S$
to a Boolean inverse semigroup is the so called Vagner--Preston map
\qcite{\Lawson}{?}
  $$
  \gamma:\S\to\I(X),
  $$
  which shows, among other things, that every inverse semigroup is a
subsemigroup of some $\I(X)$.  However $\gamma$ is never a tight
representation, even in case $\S$ is a Boolean inverse semigroup.  For
example $\gamma(0)$ is never equal to the zero of\/ $\I(X)$, namely
the empty function.  In fact this is not the only flaw presented by
$\gamma$ from the point of view of tight representations, as explained
below.

It is the main purpose of this work to introduce a canonical tight
representation
  $$
  \reg:\S \to \I(\U),
  $$
  where $\U$ is a certain space of \stress{filters}, which we call the
\stress{regular representation}. See Theorem \fcite{TheRegRepp}{6.16}.

Contrary to the Vagner--Preston representation, the regular
representation is not always faithful, but under a certain
\stress{continuity} hypothesis we are able to precisely describe when
is $\reg(s) = \reg(t)$, for a given pair of elements $s,t\in\S$. 

The issue boils down to the following situation: let $e\leq f$ be
idempotents in $E(\S)$ and suppose that there is no nonzero idempotent
$d\leq f$ such that $d\perp e$ (meaning that $de=0$).  Very roughly
speaking this means that the space between $e$ and $f$ is empty, in
which case we say that $e$ is \stress{dense} in $f$.  Notice however
that when $e\neq f$, this will never happen in a Boolean inverse
semigroup, since $d:=f\inf\neg e\neq 0$.

It turns out that when $e$ is dense in $f$ one has that
$\reg(e) = \reg(f)$, even when $e\neq f$.  In case $e$ is not
necessarily less than $f$, but $ef$ is dense in both $e$ and $f$, we
will consequently also have that $\reg(e) = \reg(ef) = \reg(f)$.

The impossibility of distinguishing between idempotents clearly has
consequences for other elements.  Suppose for example that
$s,t\in\S$ are such that $\reg(s^*s)=\reg(t^*t)$.  Suppose moreover
  that\fn{In case $\S$ is contained in some $\I(X)$, this means that
$s$ and $t$ coincide on the intersection of their domains.}
  $st^*t=ts^*s$.  Then a simple computation (see
\fcite{HowNotFaithful}{7.5}) shows that
  $\reg(s) = \reg(t)$,
  so we get another instance on non-faithfulness.

Fortunately we are able to prove in Theorem
\fcite{HowNotFaithful}{7.5} that these well understood situations are 
the only ones allowing for $\reg(s) = \reg(t)$.
  Another consequence  is that when the regular
representation is unable to separate between two elements of $\S$,
then no tight representation can possibly do it.

As already hinted upon, this result requires that $\S$ be
\stress{continuous}, as defined in \fcite{DefineContinuous}{7.1}.  To
explain what this means let us say that two elements $s,t\in\S$
\stress{essentially coincide with each other}, in symbols $s\equiv t$,
if $s^*s=t^*t$, and for every nonzero idempotent $f\leq s^*s$, there
exists a nonzero idempotent $e\leq f$, such that $se=te$.  Very
roughly this means that $s$ and $t$ coincide on a dense set, although
this idea may be made quite precise when we are speaking of
\stress{localizations} in the sense of Kumjian \cite{\Kumjian}.  See
Proposition \fcite{Localizations}{7.2}.

Recalling that when two continuous functions agree on a dense set of
their common domain they must coincide everywhere, we say that $\S$ is
\stress{continuous} if $s\equiv t$ implies that $s=t$.
  Localizations are continuous by Proposition
\fcite{Localizations}{7.2}, and so are 
  $E^*$-unitary\fn{In fact the $E^*$-unitary property, when viewed
from this point of view, reminds us of the unique continuation of
holomorphic functions: the fact that when two such functions coincide
in a small open set, they must also coincide on the largest connected
common domain of definition.  For this reason one might like to use
the expression \stress{holomorphic inverse semigroups} when referring
to the $E^*$-unitary ones.}
  inverse semigroups, as proved in \fcite{EsUnitaryCont}{7.3}.

The use of the continuity hypothesis in Theorem
\fcite{HowNotFaithful}{7.5} naturally raises the question of whether
or not this hypothesis is really needed.  To resolve this issue, in
the final section of this work we describe a general construction
which leads to a non-continuous Boolean inverse semigroup $\S$ for
which \fcite{HowNotFaithful}{7.5} does fail.

\section{Representations of semilattices}
  \label SemilatSect
  Although we are mainly interested in inverse semigroups, their
semilattice of idempotents play a particularly important role in the
ideas we shall develop.  For this reason we will set this section apart
focusing exclusively on semilattices.

\definition
  \izitem
  \zitem By a \stress{partially ordered set} we shall mean a set $X$
equipped with an order relation (i.e.~a reflexive, antisymmetric, and
transitive relation) ``$\,\leq\,$", such that $X$ contains a smallest
element, denoted $0$.
  \zitem A \stress{semilattice} is a partially ordered set $X$ 
such that for every $x,y\in X$, the set
  \ $
  \{z\in X: z\leq x,y\}
  $ \
  contains a maximum element, denoted $x\inf y$.

\medskip It is perhaps not usual to require partially ordered sets or
semilattices to contain a zero element.  However if a partially
ordered set $X$ does not contain zero one can easily embed it in
$X\mathrel{\dot\cup}\{0\}$, with the order extended from $X$ in such a
way that $0\leq x$, for all $x$.  If $X$ is a semilattice, it is
obvious that $X\mathrel{\dot\cup}\{0\}$ is also a semilattice.


\definition If $X$ is a partially ordered set we shall say that two
elements $x,y\in X$ are \stress{disjoint}, in symbols
  \ $
  x\perp y,
  $ \
  if there is no nonzero $z\in X$ such that $z\leq x,y$.  Otherwise we
shall say that $x$ and $y$ \stress{intersect}.  We shall express the
fact that $x$ and $y$ intersect by writing
  \ $
  x\its y.
  $

If $\E$ is a semilattice it is easy to see that two elements
$x,y\in\E$ intersect if and only if $x\inf y \neq 0$.

\definition 
  \label DefineRepSemilat
  Let $\E$ be a semilattice and let 
  $\B=(\B,0,1,\inf,\sup,\neg)$ 
  be a Boolean algebra.  By a \stress{representation} of $\E$ in $\B$ we
shall mean a map
  $
  \sigma:\E\to\B,
  $
  such that \izitem
  \zitem $\sigma(0) = 0$, and
  \zitem $\sigma(x\inf y) = \sigma(x)\inf\sigma(y),$ for every
$x,y\in\E$.

\medskip Recall that a Boolean algebra $\B$ is also a semilattice
  under the standard order relation given by
  $$
  \alpha\leq\beta \iff \alpha = \alpha\inf\beta
  \for \alpha,\beta\in\B.
  $$
  If $\sigma$ is a representation of the semilattice $\E$ in a Boolean
algebra $\B$ then for every $x,y\in\E$, such that $x\leq y$, one has
that $x = x\inf y$, and hence
  $$
  \sigma(x) = \sigma(x\inf y) = \sigma(x)\inf\sigma(y),
  $$
  which means that $\sigma(x)\leq \sigma(y)$.  In other words, $\sigma$
preserves the respective order relations.

\medskip
An elementary representation of any given semilattice $\E$ is obtained
as follows: let $\E^*=\E\setminus \{0\}$ and let $\P(\E^*)$ be the
Boolean algebra of all subsets of $\E^*$ under the operations of
intersection and union.  Define $\sigma:\E\to\P(\E^*)$ by setting
  $$
  \sigma(x) = \{y\in\E^*: y\leq x\}. 
  $$
  It is then easy to see that $\sigma$ is a representation of $\E$ in
$\P(\E^*)$.

Fix for the time being a representation $\sigma$ of a semilattice $\E$
in a Boolean algebra $\B$.  If $x,y\in\E$ are such that $x\leq y$, we
have already seen that $\sigma(x)\leq\sigma(y)$.
  On the other hand, if $x\perp y$, one has that
$\sigma(x)\perp\sigma(y)$, which may also be expressed in $\B$  as 
  $$
  \sigma(x)\leq \nega \sigma(y).
  $$
  More generally, if $X$ and $Y$ are finite subsets of $\E$, and one is
given an element $z\in\E$ such that $z\leq x$ for every $x\in X$, and
$z\perp y$ for every $y\in Y$, it follows that
  $$
  \sigma(z) \leq 
  \bigwedge_{x\in X} \sigma(x) \wedge 
  \bigwedge_{y\in Y} \nega{\sigma(y)}.
  \eqmark CrucialInequality
  $$
  The set of all such $z$'s will acquire an increasing importance, so we
make the following:

\definition
  \label DefineSXY
  Given finite subsets $X,Y\subseteq\E$, we shall denote by $\E^{X,Y}$
the subset of\/ $\E$ given by 
  $$
  \E^{X,Y} = \{z\in\E: z\leq x,\ \forall x\in X,\hbox{ and } z\perp
y,\ \forall y\in Y\}.
  $$

Notice that if 
  $
  \ds x_{min}=  \bigwedge_{x\in X}x,
  $
  one may replace $X$ in \lcite{\DefineSXY} by the singleton
$\{x_{min}\}$, without altering $\E^{X,Y}$.  However there does not
seem to be a similar way to replace $Y$ by a smaller set.

\definition Given any subset $F\subseteq\E$, we shall say that a
subset $Z\subseteq F$ is a \stress{cover} for $F$, if for
every nonzero $x\in F$, there exists $z\in Z$ such that $z\its x$.

The notion of covers is relevant to the introduction of the following
central concept:

\definition 
  \label DefLatTightRep  
  Let $\sigma:\E\to\B$ be a representation of the semilattice $\E$ in
the Boolean algebra $\B$.  We shall say that $\sigma$ is
\stress{tight} if for every finite subsets $X,Y\subseteq\E$, and for
every finite cover $Z$ for
  $\E^{X,Y}$,  one has that 
  $$
  \bigvee_{z\in Z}\sigma(z) \geq  \bigwedge_{x\in X} \sigma(x) \wedge
\bigwedge_{y\in Y} \nega{\sigma(y)}.
  $$

Notice that the reverse inequality ``$\leq$" always holds by
\lcite{\CrucialInequality}.  Thus, when $\sigma$ is tight, we actually
get an equality above.  We should also remark that in the absence of
any finite cover $Z$, as above, every representation is considered to
be tight by default.

In certain cases the verification of tightness may be greatly simplified:

\state Proposition
  \label AltLatTightRep
  Let $\sigma$ be a representation of the semilattice $\E$ in the
Boolean algebra $\B$, such that either 
  \izitem 
  \zitem $\E$ contains a finite set $X$ such that $\bigvee_{x\in X}
\sigma(x) = 1$, or
  \zitem $\E$ does not admit any finite cover.
  \medskip\noindent
  Then $\sigma$ is tight if and only if for every nonzero $x\in\E$ and
for every finite
  cover $Z$ for the \stress{interval}
  $$
  [0, x]:= \{z\in\E: z\leq x\}, 
  $$
  one has that
  $
  \bigvee_{z\in Z}\sigma(z) \geq \sigma(x).
  $

\proof
  See \scite{\actions}{10.8}.
  \proofend

Whenever $z\in[0, x]$, notice that $\sigma(z) \leq \sigma(x)$, so the
last inequality in the statement of the result above is in fact
equivalent to $\bigvee_{z\in Z}\sigma(z) =\sigma(x)$.

The representation of $\E$ in $\P(\E^*)$ described above is not
necessarily tight.  In fact, if $\E$ consists of three distinct
  elements, say $\E=\{0,y,1\}$, with the order relation such that
$0\leq y\leq1$, set
$X=\{1\}$ and $Y=\{y\}$.  Then $\E^{X,Y}=\{0\}$, so the empty set $Z$
is a cover for $\E^{X,Y}$.  However
  $$
  \bigvee_{z\in Z}\sigma(z) = 
  \emptyset \neq 
  \{1\} = 
  \bigwedge_{x\in X} \sigma(x) \wedge \bigwedge_{y\in Y}
\nega{\sigma(y)}.
  $$

Not all semilattices admit tight injective representations.  In order to
study this issue in detail it is convenient to introduce the following:

\definition 
  \label DefineDense
  Let $\E$ be a semilattice and let $x,y\in\E$ be such that $y\leq x$.
We shall say that $y$ is \stress{dense} in $x$ if there is no nonzero
$z\in\E$ such that $z\leq x$ and $z\perp y$.  Equivalently, if
$\E^{\{x\},\{y\}}=\{0\}$.

Obviously each $x\in\E$ is dense in itself but it is conceivable that
some $y\neq x$ is dense in $x$.  For a concrete example notice that in
the semilattice $\E=\{0,y,1\}$ above one has that $y$ is dense in $1$.

In the general case, whenever $y$ is dense in $x$ we have that
  $
  \E^{\{x\},\{y\}}=\{0\},
  $
  and hence the empty set is a cover for 
  $\E^{\{x\},\{y\}}$.  Therefore for every tight representation $\sigma$
of $\E$ one has that
  $$
  0=\sigma(x) \wedge \nega{\sigma(y)},
  $$
  which means that $\sigma(x) \leq \sigma(y)$.  Since the opposite
inequality also holds, we have that $\sigma(x)=\sigma(y)$.  Thus no
tight representation of $\E$ can possibly separate $x$ and $y$.
For future reference we record this conclusion in the next:

\state Proposition
  \label ThouShallNotSeparate
  If $y\leq x$ are elements in the semilattice $\E$, such that $y$ is
dense in $x$, then $\sigma(y)=\sigma(x)$ for every tight
representation $\sigma$ of\/ $\E$.

 %
 %
 %
 %

When $\E$ happens to be a Boolean algebra there is a very
elementary characterization of tight representations:

\state  Proposition
  \label BooleToBoole
  Suppose that $\E$ is a semilattice admitting the structure of a
Boolean algebra which induces the same order relation as that of $\E$,
and let 
  \ $
  \sigma:\E\to\B 
  $ \
  be a representation of $\E$ in some Boolean algebra $\B$.  Then
$\sigma$ is tight if and only if it is a Boolean algebra homomorphism.

\proof
  Supposing that $\sigma$ is tight, notice that $\{1\}$ is a cover for
$\E^{\emptyset,\{0\}}$, so
  $$
  \sigma(1) = \neg \sigma(0) = \neg 0 = 1.
  $$
  Given $x\in\E$ notice that $\{\neg x\}$ is a cover for
$\E^{\emptyset,\{x\}}$, therefore 
  $$
  \sigma(\neg x) = \neg\sigma(x).
  $$
Since 
  $
  x\sup y =   \neg(\neg x\inf \neg y),
  $
  for all $x,y\in\E$,  we may easily prove that $\sigma(x\sup
y)=\sigma(x)\sup \sigma(y)$.  Thus $\sigma$ is a Boolean algebra
homomorphism, as required.

In order to prove the converse implication let $X,Y\subseteq\E$ be
finite sets and let $Z$ be a finite cover for $\E^{X,Y}$.  Let
  $$
  z_0 = \bigvee_{z\in Z} z, \quad
  x_0 = \bigwedge_{x\in X} x \and
  \bar y_0 = \bigwedge_{y\in Y} \neg y.
  $$
  It is obvious that $z_0\leq x_0\inf \bar y_0$, and we claim that in
fact $z_0= x_0\inf \bar y_0$.  We will prove it by checking that
  $$
  \neg z_0\inf x_0\inf \bar y_0 = 0.  
  $$
  Let $u = \neg z_0\inf x_0\inf \bar y_0$, and notice that the fact that
$u\leq x_0\inf \bar y_0$ implies that $u\in\E^{X,Y}$. Arguing by
contradiction, and hence supposing that $u$ is nonzero, we deduce that
$u\its z$, for some $z\in Z$, but this contradicts the fact that $u\leq
\neg z_0$.  This proves our claim so, assuming that $\sigma$ is a
Boolean algebra homomorphism, we have
  $$
  \bigvee_{z\in Z} \sigma(z) =
  \sigma\Big( \bigvee_{z\in Z} z\Big) =
  \sigma(z_0) =
  \sigma(x_0\inf \bar y_0) =
  \bigwedge_{x\in X} \sigma(x) \inf 
  \bigwedge_{y\in Y} \neg \sigma(y),
  $$
  showing that $\sigma$ is tight.
  \proofend

We shall have a lot more to say about tight representations in the
following sections.

\section{Filters}
  A fundamental tool for the study of tight representations of
semilattices is the notion of filters, which we briefly introduce in
this section.

\definition
  \label DefineFilter
  Let $X$ be any partially ordered set with minimum element $0$.
  A \stress{filter} in $X$ is a nonempty
subset $\xi\subseteq X$, such that
  \izitem
  \zitem $0\notin\xi$,
  \zitem if $x\in\xi$ and $y\geq x$, then $y\in\xi$,
  \zitem if $x,y\in\xi$, there exists $z\in\xi$, such that $x,y\geq z$.
  \medskip \noindent An \stress{ultrafilter} is a filter which is not
properly contained in any filter. 

Given a partially ordered set $X$ and any nonzero element $x\in X$ it is
elementary to prove that
  $$
  \xi=\{y\in X: y\geq x\}
  $$
  is a filter containing $x$.  By Zorn's Lemma there exists an
ultrafilter containing $\xi$,   thus every nonzero element in $X$ belongs
to some ultrafilter.

When $\E$ is a semilattice, given the existence of $x\inf y$ for every
$x,y\in\E$, condition \lcite{\DefineFilter.iii} may be replaced by 
  $$
  x,y\in\xi \  \imply  \ x\inf y\in\xi.
  \eqmark NewCondForFilInSL
  $$

The following is an important fact about filters in semilattices which
also benefits from the existence of $x\inf y$.

\state Lemma
  \label UltrafilterCriterium
  Let $\E$ be a semilattice and let $\xi$ be a filter in $\E$.  Then
$\xi$ is an ultrafilter if and only if $\xi$ contains every element
$y\in\E$ such that $y\its x$ for every $x\in\xi$.

\proof
  In order to prove the ``if" part let $\eta$ be a filter such that
$\xi\subseteq\eta$.  Given $y\in\eta$ one has that for every
$x\in\xi$, both $y$ and $x$ lie in $\eta$, and hence 
\lcite{\NewCondForFilInSL} implies that
  $y\inf x\in\eta$, so $y\inf x\neq 0$, and hence $y\its x$.  By
hypothesis $y\in\xi$, proving that $\eta=\xi$, and hence that $\xi$ is
an ultrafilter.

Conversely let $\xi$ be an ultrafilter and suppose that $y\in\E$ is
such that $y\its x$, for every $x\in\xi$.  Defining
  $$
  \eta = \{ u\in\E: u \geq y\inf x, \hbox{ for some } x\in \xi\},
  $$
  we claim that $\eta$ is a filter.
  By hypothesis $0\notin\eta$. 
  Also if $u_1,u_2\in\eta$, choose for every $i=1,2$ some $x_i\in\xi$
such that $u_i\geq y\inf x_i$.  Then
  $$
  u_1\inf u_2 \geq (y\inf x_1)\inf(y\inf x_2) = y\inf (x_1\inf x_2),
  $$
  so $u\in\eta$.  Given that \lcite{\DefineFilter.ii} is obvious we
see that $\eta$ is indeed a filter, as claimed.  Noticing that
$\xi\subseteq\eta$ we have that $\eta=\xi$, because $\xi$ is an
ultrafilter.  Since $y\in\eta$, we deduce that $y\in\xi$.
  \proofend

\section{Characters}
  \label SpecSection
  We fix, throughout this section, a semilattice $\E$, always assumed
to have a smallest element $0$.  The study of representations of $\E$
in the most elementary Boolean algebra of all, namely $\{0,1\}$, leads
us to the following well known  important concept.

\definition
  A \stress{character} of $\E$ is a nonzero  representation of $\E$ in
the Boolean algebra $\{0, 1\}$.
The set
of all characters will be  denoted by $\spec$.

Some authors use the term \stress{semicharacter} referring to maps 
$\phi:\E\to\{0,1\}$ satisfying \lcite{\DefineRepSemilat.ii}.
Thus, a character is nothing but a semicharacter which vanishes at
$0$.  Perhaps the widespread use of the term \stress{semicharacter} is
motivated by the fact that it shares prefix with the term
\stress{semilattice}.  If this is really the case then our choice of
the term \stress{character} may not be such a good idea but alas, we
cannot think of a better term.

Temporarily denoting by $\widetilde\E$ the
set of \stress{all} representations of $\E$ in $\{0,1\}$, including
the identically zero representation, is is easy to see that
$\widetilde\E$ is a closed subspace of the compact product space
$\{0,1\}^\E$, hence $\widetilde\E$ is compact.  Since $\spec$ is
obtained by removing the identically zero representation from
$\widetilde\E$, we have that $\spec$ is locally compact.

Given a character $\phi$, observe that 
  $$
  \xi_\phi  = \{x\in\E: \phi(x) = 1\},
  \eqmark XiPhi
  $$
  is a filter in $\E$ (it is nonempty because $\phi$ is assumed not to
be identically zero).
   Conversely, given a filter $\xi$, define for every $x\in\E$,
  $$
  \phi_\xi(x) = \left\{\matrix{
  1, & \hbox{ if } x\in\xi,\hfill\cr
  \pilar{12pt}
  0, & \hbox{ otherwise.}}\right.
  $$
  It is then easy to see that $\phi_\xi$ is a character.
Therefore we see that \lcite{\XiPhi} gives a one-to-one correspondence
between $\spec$ and the set of all filters.

\state  Proposition
  \label UltraThenTight
  If $\xi$ is an ultrafilter then $\phi_\xi$ is a tight representation
of\/ $\E$ in $\{0,1\}$.

  \proof
  Let $X,Y\subset \E$ be finite subsets and let $Z$ be a cover for
$\E^{X,Y}$.  In order to prove that
  $$
  \bigvee_{z\in Z}\phi(z) \geq 
  \prod_{x\in X} \phi(x) \prod_{y\in Y} (1-\phi(y)),
  $$
  it is enough to show that if the right-hand side equals 1, then so
do the left-hand side.  This is to say that if $x\in\xi$ for every
$x\in X$, and $y\notin\xi$ for every $y\in Y$, then there is some
$z\in Z$, such that $z\in\xi$. 

By \lcite{\UltrafilterCriterium}, for each $y\in Y$ there exists some
$x_y\in\xi$ such that $y\perp x_y$.
Supposing by contradiction that
$Z\cap\xi = \emptyset$, then for every $z\in Z$ there exists, again by
\lcite{\UltrafilterCriterium}, some $x_z\in\xi$, such that $z\perp
x_z$.  Set 
  $$\hbox{$
  w =
  \bigwedge\limits_{x\in X} x \wedge
  \bigwedge\limits_{y\in Y} x_y \wedge
  \bigwedge\limits_{z\in Z} x_z.
  $}$$
  Since $w\in\xi$ we have that $w\neq 0$.  Obviously $w\leq x$ for
every $x\in X$, and $w\perp y$ for every $y\in Y$, and hence
$w\in\E^{X,Y}$.  Since $Z$ is a cover there exists some $z_1\in Z$
such that $w\its z_1$.  However, since $w\leq x_{z_1}\perp z_1$, we
have that $w\perp z_1$, a contradiction.
  \proofend

\definition
  \label AltSpectrums
  We shall denote by $\speci$ the set of all characters
$\phi\in\spec$ such that $\xi_\phi$ is an ultrafilter.  Also we will
denote by $\spect$ the set of all tight characters.

Employing the terminology just introduced we may rephrase
\lcite{\UltraThenTight} by saying that $\speci\subseteq\spect$.
The following main result further describes the relationship
between $\speci$ and $\spect$.

\state Theorem
  \label ClosureOfUltrafilters
  Let $\E$ be a semilattice with smallest element $0$, and let
$\speci$ and $\spect$ be as defined in \lcite{\AltSpectrums}.  Then
the closure of $\speci$ in $\specz$ coincides with $\spect$.

\proof Since the condition for any given  $\phi$ in $\specz$ to belong to
$\spect$ is given by equations 
it is easy to prove that $\spect$ is closed within $\specz$,
and since $\speci\c\spect$ by \lcite{\UltraThenTight}, we
deduce that
  $$
  \overline{\speci}\c\spect.
  $$
  To prove the reverse inclusion let us be given $\phi\in\spect$.  We
must therefore show that $\phi$ can be arbitrarily approximated by
elements from $\speci$.  Let $U$ be a neighborhood of $\phi$
within $\specz$.  By definition of the product topology, $U$ contains
a neighborhood of $\phi$ of the form
  $$
  V = V_{X,Y} = 
  \{\psi\in\specz: 
  \psi(x) = 1, \hbox{ for all } x\in X, \hbox{ and }
  \psi(y) = 0, \hbox{ for all } y\in Y\},
  $$
  where $X$ and $Y$ are finite subsets of $\E$.  
We next claim that $\E^{X,Y}\neq \{0\}$.  In order to prove this
suppose the contrary, and hence $Z=\emptyset$ is a cover for
$\E^{X,Y}$.  Since $\phi$ is tight we conclude that 
  $$
  0 = \bigvee_{z\in Z}\phi(z) =
  \prod_{x\in X} \phi(x) \prod_{y\in Y} (1-\phi(y)).
  $$
  However, since $\phi$ is supposed to be in $V$, we have that 
$\phi(x)=1$ for all $x\in X$, and $\phi(y)=0$ for all $y\in Y$, which
means that the right-hand side of the expression displayed above
equals 1.  This is a contradiction and hence our claim is proved.

We are therefore allowed to choose a nonzero $z\in \E^{X,Y}$, and
further to pick
an ultrafilter $\xi$ such that $z\in\xi$.  Observe that
$\phi_\xi\in\speci$, and the proof will be concluded once we
show that $\phi_\xi\in U$.

For every $x\in X$ and $y\in Y$, we have that $z\leq x$ and $z\perp
y$, hence $x\in\xi$ and $y\notin \xi$.  This entails $\phi_\xi(x)=1$
and $\phi_\xi(y)=0$, so $\phi_\xi\in V\c U$, as required.
  \proofend

In the correspondence between $\spec$ and the set of all filters given by
\lcite{\XiPhi}, we know that elements of $\speci$ correspond to
ultrafilters (by definition of the former).  Given the importance of
the notion of tight characters, highlighted by
\lcite{\ClosureOfUltrafilters}, it is sensible to make the following:

\definition
  A filter $\xi$ in $\E$ is said to be \stress{tight} if
$\phi_\xi$ is a tight character, that is, if $\phi_\xi\in\spect$.

By \lcite{\UltraThenTight} we see that every ultrafilter is tight.


  %

\section{Boolean inverse semigroups}
  Recall that a semigroup $\S$ is said to be an
\stress{inverse semigroup} if for every $s\in\S$, there exists a
unique $s^*\in\S$ such that
  $$
  ss^*s=s
  \and
  s^*ss^*=s^*.
  $$
  See \cite{\Lawson} for a detailed study of inverse semigroups.
  It is well known that the correspondence $s\mapsto s^*$ is then an
involutive anti-homomorphism.

  One usually denotes by $E(\S)$ the set of all idempotent elements of
$\S$, such as $s^*s$, for every $s\in\S$.  It is not hard to show that
$E(\S)$ is a semilattice under the order
  $$
  e\leq f\iff e = ef
  \for e,f\in E(\S).
  $$
  In particular one has that $e\inf f=ef$, for all $e,f\in E(\S)$.

  A \stress{zero} element of a semigroup $\S$ is by definition an
element, usually denoted $0$, such that
  $$
  0s = s0 = 0 \for s\in\S.
  \eqmark ZeroElement
  $$
  Any semigroup $\S$ can be readily embedded in a semigroup with zero
by simply adding an extra element, denoted $0$, and extending the
multiplication operation of $\S$ by means of \lcite{\ZeroElement}.  If $\S$
happens to be 
an inverse semigroup it is easy to show that $\S\cup\{0\}$ is also an
inverse semigroup, with $0^*=0$.


Given that semigroups with zero are often difficult to
handle, one may wonder why in the world would anyone want to insert a
zero in an otherwise well behaved semigroup.  
  Rather than shy away from inverse semigroups with zero, we will
assume that \stress{all} of them contain a zero element, not least
because we want to keep a close eye on this exceptional element.

In addition to the standard order on $E(\S)$ described above it is
important to consider a certain order relation on $\S$.

\definition
  \label DefineOrderInISG
  Given an inverse semigroup $\S$, and 
  given $s,t\in\S$, we will say that $s\leq t$ \ if any one of the
following equivalent relations hold:
  \izitem
  \zitem $ts^*s = s$,
  \zitem $ss^*t = s$,
  \zitem there exists $e\in E(\S)$ such that $te=s$,
  \zitem there exists $e\in E(\S)$ such that $et=s$.

See \qcite{\Lawson}{???} for a proof of the fact that these conditions are in
fact equivalent.

\bigskip

  We have already mentioned that every Boolean algebra $\B$ is a partially
ordered set with the order defined by
  $x\leq y$ if $ x\inf y = x$.
  This order in fact encodes all of the Boolean algebra structure of
$\B$ since
  $$
  \eqmatrix{
  \hfill 0 &=& \min\B, \hfill\cr
  \hfill 1 &=& \max\B, \hfill\cr
  x\inf y &=& \min\{x,y\}, \hfill\cr
  x\sup y &=& \max\{x,y\}, \hfill\cr
  \hfill \neg x  &=& \max\{y: x\inf y = 0\}.}  
  $$
  Thus, if $\B_1$ and $\B_2$ are Boolean algebras and
$\phi:\B_1\to\B_2$ is an order-isomorphism, i.e., a bijective map such
that
  $$
  x\leq y\iff \phi(x)\leq\phi(y)
  \for x,y\in\B_1,
  $$
  then $\phi$ is in fact a Boolean algebra isomorphism.

\definition
  Given any partially ordered set $X$ we will say that $X$
\stress{is} a Boolean algebra if $X$ is order-isomorphic to a
(necessarily unique) Boolean algebra.

Suppose that $\S$ is an inverse semigroup whose semilattice of
idempotents $E(\S)$ is a Boolean algebra in the above sense.  In
particular $E(\S)$ must contain a smallest element 0, and a biggest
element 1.  For every $e\in E(\S)$ one then has that $1e=e$.  It
follows that for every $s\in\S$,
  $$
  1s = 1ss^*s = ss^*s = s,
  $$
  and similarly $s1=s$.  So we see that 1 is a multiplicative unit for $\S$.

On the other hand $0e=e0=0$, for every idempotent $e$, but this does
not necessarily imply \lcite{\ZeroElement}, a counter-example being
that of any group $G$ with more than one element.  In this case the
smallest element of $E(G)=\{1\}$ \ is $1$, so $0=1$.  However it is
definitely not true that $0s=s0=0$, for every $s\in G$.


\definition
  A \stress{Boolean inverse semigroup} is an inverse semigroup
$\BS$ whose lattice of idempotents $E(\BS)$ is a Boolean algebra, and
such that
  $$
  0s=s0=0
  \for s\in\BS,
  $$
  where $0$ denotes the smallest element of $E(\BS)$.

If $\BS$ is a Boolean inverse semigroup we will freely use Boolean
algebra language when referring to $E(\BS)$ with the understanding that
it relates to the unique Boolean algebra structure compatible with the
order structure of $E(\BS)$.  However we will refrain from using the
\stress{meet} operator ``$\inf$" since it is conveniently substituted
by the multiplication operation on $E(\BS)$.  

Examples of Boolean inverse semigroups are quite common.  If $X$ is
any set and $\I(X)$ is the set of all partially defined bijections
between subsets of $X$ then it is well known that $\I(X)$ is an
inverse semigroup under composition.  The semilattice of idempotents
of $\I(X)$ is identical to $\P(X)$, the Boolean algebra of all subsets
of $X$, and hence $\I(X)$ is a Boolean inverse semigroup.

\section{Representations of inverse semigroups}
  \label TightrRpSection
  Throughout this section we fix an inverse semigroup $\S$ with zero.
If $\B$ is a Boolean inverse semigroup and
  $$
  \sigma:\S\to\B
  $$
  is a semigroup homomorphism, observe that $\sigma(E(\S))\subseteq
E(\B)$, so the restriction of $\sigma$ to $E(\S)$ is a map from a
semilattice into a Boolean algebra.  It therefore makes sense to ask
whether or not it is a tight representation.

\definition 
  \label DefISGTightRep
  Let $\BS$ be a Boolean inverse semigroup.  A semigroup homomorphism
  $\sigma: \S \to \BS$
  is said to be a \stress{tight representation} if the restriction of
$\sigma$ to $E(\S)$ is a tight representation of $E(\S)$ in $E(\BS)$,
in the sense of \lcite{\DefLatTightRep}.

Notice that if $\sigma$ is a tight representation in the above sense
then \lcite{\DefineRepSemilat.i} applies so it is understood that
$\sigma(0)=0$.  We also remind the reader that a semigroup
homomorphism $\sigma$ between inverse semigroups necessarily satisfies
  $
  \sigma(s^*) = \sigma(s)^*.
  $

The following is an obvious consequence of \lcite{\BooleToBoole}.

\state Proposition
  \label BisgToBisg
  Let $\S$ and $\B$ be Boolean inverse semigroups and let
  $$
  \sigma:\S\to\B 
  $$
  be a semigroup
homomorphism.  Then  $\sigma$ is a tight representation if and only if
the restriction of $\sigma$ to $E(\S)$ is a homomorphism in the
category of Boolean algebras.

Among the better known examples of a semigroup homomorphism from an
inverse semigroup $\S$ to a Boolean inverse semigroup is the
Vagner--Preston representation \qcite{\Lawson}{?}, so it is interesting to
ask whether or not it is a tight representation.  In order to fix
notation let us briefly describe it.  For every idempotent $e\in
E(\S)$, let
  $$
  D_e = \{t\in\S: tt^* \leq e\},
  $$
  and for $s\in\S$ consider the map
  $
  \gamma(s):D_{s^*s} \to D_{ss^*},
  $
  given by $\gamma(s)(t)=st$.  Then each $\gamma(s)$ is a bijective map
and hence $\gamma$ gives a map
  $$
  \gamma:\S \to \I(\S),
  $$
  which is well known to be a semigroup monomorphism.  The
Vagner--Preston Theorem asserts that every inverse semigroup may be
realized inside some $\I(X)$ and $\gamma$ provides just that
realization.

Supposing, as we are, that $\S$ contains a zero element, notice that
$D_0$ is the singleton $\{0\}$, while $\gamma(0)$ is the identity map
on $D_0$, so that $\gamma(0)$ is \stress{not} the zero element of
$\I(\S)$, the latter being the empty function.  This violates
\lcite{\DefineRepSemilat.i} and hence $\gamma$ is \stress{not} a tight
representation.  One could remedy this by removing zero from every
$D_e$, but it would still not give us a tight representation.
  To see this consider, for example, the following Boolean algebra
viewed as a semilattice, and hence as an inverse semigroup:
  $$
  \S = \{0,1\}\times\{0,1\}.
  $$
  Since all elements of $\S$ are idempotent, the range of $\gamma$ is
contained in $E(\I(\S))=\P(\S)$.   Removing zero as suggested
above, $\gamma$ becomes  the map
  $$
  \eqmatrix{
  (0,0) & \to & \emptyset \cr
  (1,0) & \to & \{(1,0)\}\cr
  (0,1) & \to & \{(0,1)\}\cr
  (1,1) & \to & \{(1,0),(0,1),(1,1)\}.
  }
  $$
  Observe that if $X=\{(1,1)\}$, and $Y=\emptyset$, then
$E(\S)^{X,Y} = E(\S) = \S$, and hence 
$Z=\{(1,0),(0,1)\}$ is a cover for $E(\S)^{X,Y}$.  However
  $$
  \bigvee_{z\in Z}\gamma(z) = \gamma(1,0)\cup \gamma(0,1) =
\{(1,0),(0,1)\},
  $$
  while 
  $$
  \bigwedge_{x\in X} \gamma(x) \wedge
  \bigwedge_{y\in Y} \nega{\gamma(y)} =
  \gamma(1,1) =  \{(1,0),(0,1),(1,0)\}.
  $$

  The purpose of this section is to exhibit a canonical tight
representation of $\S$.  Filters will again be crucial in 
achieving this.   Whenever we speak of filters in $\S$ it will be with
respect to the standard order
relation on $\S$ given by \lcite{\DefineOrderInISG},

  If $\xi$ is a filter in $\S$ and $e\in E(\S)$, let
  $$
  e\xi = \{et: t\in\xi\}.
  $$
  We will now turn our attention to filters $\xi$ such that
$e\xi\subseteq\xi$.

\state Lemma
  \label CoolLemma
  Given a filter $\xi$ suppose that $es\in\xi$, for some $s\in\S$ and
$e\in E(\S)$.  Then $e\xi\subseteq\xi$.

\proof
  Given $t\in\xi$ observe that by \lcite{\DefineFilter.iii} there
exists $r\in\xi$ such that $es,t\geq r$.  Therefore $r= esr^*r$, so
  $$
  e t \geq e r = e (esr^*r) = esr^*r = r\in\xi,
  $$
  so
  $
  e t\in\xi.
  $
  \proofend

\state Corollary
  \label EEtainEta
  If $\xi$ and $\eta$ are filters such that $\xi\subseteq\eta$, and
$e$ is an idempotent such that $e\xi\subseteq\xi$, then
$e\eta\subseteq\eta$.

\proof
  Given any $s\in\xi$ we have that $es\in e\xi\subseteq\xi\subseteq\eta$, and
hence $e\eta\subseteq\eta$ by \lcite{\CoolLemma}.

\state Corollary
  If $\xi$ is a filter and $s\in\xi$, then $ss^*\xi\subseteq\xi$.

\proof
  Since $ss^*s = s \in \xi$, the result follows from
\lcite{\CoolLemma}.
  \proofend

For ultrafilters there is another important condition which implies
the same conclusion as \lcite{\CoolLemma}:

\state Lemma
  \label UltraCool
  Let $\xi$ be an ultrafilter and let $e\in E(\S)$ be such that $et\neq 0$,
for all $t\in\xi$.  Then $e\xi\subseteq\xi$. 

\proof
  Let 
  $$
  \eta = \{u\in\S: u\geq et, \hbox{ \sl for some } t\in\xi\}.
  $$
  Observe that $\eta$ is a filter, since $0\notin\eta$, by hypothesis,
and \lcite{\DefineFilter.ii-iii} are of easy verification.  For every
$t\in\xi$ one has that $t\geq et$, and hence $t\in\eta$.  Thus
$\xi\subseteq\eta$, and since $\xi$ is an ultrafilter, we deduce that
$\xi=\eta$.  For $t\in\xi$, it is obvious that $et\in\eta$.  This says
that $e\xi\subseteq\eta=\xi$.
  \proofend


\definition
  \label DefineFilterSpc
  We will denote by $\U$ the set of all ultrafilters in $\S$, and for
each idempotent $e\in E(\S)$ we will denote by $\U_e$ the set of
all ultrafilters $\xi$ such that $e\xi\subseteq\xi$.

The following result describes how do the $\U_e$ behave under
intersections.

\state Lemma 
  \label FilterIntersection
  Let $e$ and $f$ be idempotents in $E(\S)$.  Then \
  $
  \U_e \cap \U_f = \U_{ef}.
  $

\proof
  If $\xi\in\U_e \cap \U_f$ then 
  $$
  ef\xi=e(f\xi)\subseteq e\xi \subseteq \xi,
  $$
  so $\xi\in\U_{ef}$.  Conversely, if $\xi$ is in $\U_{ef}$, pick any
$t$ in $\xi$.  Then $eft\in ef\xi\subseteq\xi$, so we deduce from
\lcite{\CoolLemma} that $e\xi\subseteq\xi$, and hence $\xi\in\U_e$.
Similarly $fet\in \xi$, so $\xi\in\U_f$, as desired.
  \proofend

Notice that if $e=0$, then $\U_e=\emptyset$.  Thus the above result
implies that $\U_e$ and $\U_f$ are disjoint when $ef=0$.

Let us now take some time to discuss when is $\U_e=\U_f$, for
idempotents $e$ and $f$.

\state Proposition
  \label SameUs
  Let $e$ and $f$ be idempotents in $E(\S)$.  Then $\U_e=\U_f$ if and
only if \ $ef$ \ is dense \lcite{Defi\-nition \DefineDense} in both
$e$ and $f$.  In this case, for every tight representation $\sigma$
of\/ $\S$, one has that $\sigma(e)=\sigma(f)$.

\proof
  Let us first prove the only if part.  We begin by treating the
special case in which $e\leq f$.  Thus, assuming that $\U_e=\U_f$,
we must prove that $e$ is dense in $f$. 

Arguing by contradiction, let $d$ be a nonzero idempotent such that
$d\perp e$, and $d\leq f$.  Choose an ultrafilter $\xi$ such that
$d\in\xi$ and observe that 
  $$
  fd = d\in\xi,
  $$
  so that $\xi\in\U_f$, by \lcite{\CoolLemma}.  By assumption we have
that $\xi\in\U_e$ and hence $e\xi\subseteq\xi$.  In particular
  $$
  0=ed \in e\xi \subseteq\xi,
  $$
  which is a contradiction.  This proves that $e$ is dense in $f$.

Without the assumption that $e\leq f$, but still supposing that
$\U_e=\U_f$, observe that by \lcite{\FilterIntersection} we have
  $$
  \U_e = \U_{ef} = \U_f.
  $$
  By the first part of the proof we then deduce that $ef$ is dense in
both $e$ and $f$, as required. 

Conversely, suppose that $ef$ is dense in $e$ and $f$.  In order to
conclude the proof it is obviously enough to prove that
$\U_e=\U_{ef}$, and $\U_f=\U_{ef}$, while by symmetry it suffices to
prove only the first assertion.  Observing that $\U_{ef}\subseteq\U_e$
by \lcite{\FilterIntersection},  we must only prove that
$\U_e\subseteq\U_{ef}$.

For this let $\xi\in\U_e$.  Given any $t\in\xi$ we claim that
$eft\neq0$.  To prove it suppose otherwise so that $eftt^* = 0$, for
some $t\in\xi$.  This says that $ett^*\perp ef$, and clearly
$ett^*\leq e$.  Since $ef$ is dense in $e$ we deduce that $ett^*=0$,
hence
  $$
  0 = ett^*t = et \in e\xi\subseteq\xi,
  $$
  which is a contradiction.  This shows that $eft\neq0$, for every
$t\in\xi$.  By \lcite{\UltraCool} it then follows that
$\xi\in\U_{ef}$, as desired.

Finally, given a tight representation of $\S$ we have that the
restriction of $\sigma$ to $E(\S)$ is a tight representation of
$E(\S)$ in the sense of \lcite{\DefLatTightRep}.  Hence we have by
\lcite{\ThouShallNotSeparate} that
  $$
  \sigma(e) = \sigma(ef) = \sigma(f),
  $$
  proving the last part of the statement.
  \proofend

With the next result we shall start to study certain functions on the
set of filters, in preparation for introducing the regular 
representation. 

\state Proposition
  \label MovingFilters
  Given $s\in\S$ and a filter $\xi$ such that $s^*s\xi\subseteq \xi$,
let
  $$
  \reg_s(\xi) = 
  \{u\in\S: u\geq st, \hbox{ \sl for some } t\in\xi\}.
  $$
  Then \izitem
  \zitem $\reg_s(\xi)$ is a filter,
  \zitem $ss^*\reg_s(\xi)\subseteq\reg_s(\xi)$,
  \zitem $s\xi\subseteq\reg_s(\xi)$.

\proof 
  With respect to the last assertion let $t\in\xi$, and put $u=st$.
Then obviously $u\geq st$, so $u\in\reg_s(\xi)$.  In order to prove
(i) assume by contradiction that $0\in \reg_s(\xi)$.  Then $st=0$,
for some $t\in\xi$, and hence
  $$
  0 = s^*st \in s^*s\xi\subseteq \xi,
  $$
  a contradiction, proving that $0\notin \reg_s(\xi)$.  If
$u_1,u_2\in\reg_s(\xi)$, choose for $i=1,2$, some $t_i\in\xi$ such
that $u_i\geq st_i$.  Pick $t\in\xi$ such that $t_1,t_2\geq t$, and
set $u=st$.  By (iii) one has that $u\in\reg_s(\xi)$ and we have
  $$
  u_i \geq st_i \geq st = u,
  $$
  proving \lcite{\DefineFilter.iii}.  Since \lcite{\DefineFilter.ii}
is obvious we have concluded the proof that $\reg_s(\xi)$ is a filter.

In order to prove (ii) let $u\in\reg_s(\xi)$, and pick $t\in\xi$
such that $u\geq st$.  Then
  $$
  ss^* u \geq ss^* st = st,
  $$
  so $ss^* u \in\reg_s(\xi)$.
  \proofend

Given a filter $\xi$ such that $s^*s\xi\subseteq\xi$, we have seen
above that $ss^*\reg_s(\xi)\subseteq\reg_s(\xi)$, so it makes
sense to speak of $\reg_{s^*}\big(\reg_s(\xi)\big)$.

\state Proposition
  \label VaiEVolta
  Let $\xi$ be a filter such that $s^*s\xi\subseteq \xi$.  Then
$\reg_{s^*}\big(\reg_s(\xi)\big) = \xi$.

\proof
  If $v\in \reg_{s^*}\big(\reg_s(\xi)\big)$, there exists $u\in
\reg_s(\xi)$ such that $v\geq s^*u$.  In turn there exists $t\in\xi$
such that $u\geq st$, so
  $$
  v \geq s^*u\geq s^*st \in s^*s\xi\subseteq\xi,
  $$
  and hence $v\in\xi$.  Conversely let $t\in\xi$, then
$st\in\reg_s(\xi)$ by \lcite{\MovingFilters.iii}, and by the same
token $s^*st\in\reg_{s^*}\big(\reg_s(\xi)\big)$.  Since $t\geq
s^*st$, we have that $t\in\reg_{s^*}\big(\reg_s(\xi)\big)$.
  \proofend

We shall next prove that $\reg_s$ preserves ultrafilters.  For this
recall that for $e\in E(\S)$, we denote by $\U_e$ the set of all
ultrafilters $\xi$ such that $e\xi\subseteq \xi$.

\state Proposition 
  \label MovingUltraFilters
  If  $s\in\S$ and  $\xi\in\U_{s^*s}$, then $\reg_s(\xi)\in\U_{ss^*}$.
  
\proof In order to prove that $\reg_s(\xi)$ is an ultrafilter, suppose
that $\reg_s(\xi)\subseteq\eta$, for some filter $\eta$.  We must show
that $\reg_s(\xi)=\eta$.
  By \lcite{\MovingFilters.ii} we have that
$ss^*\reg_s(\xi)\subseteq\reg_s(\xi)$, so we may use
\lcite{\EEtainEta} to conclude that $ss^*\eta\subseteq\eta$.  Thus
$\reg_{s^*}(\eta)$ is a filter by \lcite{\MovingFilters.i}.  Using
\lcite{\VaiEVolta} we have
  $$
  \xi = \reg_{s^*}(\reg_s(\xi))\subseteq\reg_{s^*}(\eta),
  $$
  so $\xi =\reg_{s^*}(\eta)$, by maximality.  This implies that 
  $$
  \reg_s(\xi) =\reg_s\big(\reg_{s^*}(\eta)\big) =\eta.
  $$
  To prove that $\reg_s(\xi)\in\U_{ss^*}$ it then suffices to show
that $ss^*\reg_s(\xi)\subseteq\reg_s(\xi)$, which is nothing but
\lcite{\MovingFilters.ii}.
  \proofend
  
The following is a useful characterization of $\reg_s(\xi)$ when
$\xi$ is an ultrafilter.

\state Proposition 
  \label UsefulCharac
  Let $s\in\S$ and let $\xi\in\U_{s^*s}$.  Then $\reg_s(\xi)$ is the
unique filter containing $s\xi$.

\proof
  By \lcite{\MovingFilters.iii} we have that $\reg_s(\xi)$ does
indeed contain $s\xi$.  So let $\eta$ be another filter such that
$s\xi\subseteq\eta$.   We must prove that $\eta=\reg_s(\xi)$.
 Given any $t \in\xi$ we have that
  $$
  ss^*st =
  st \in s\xi \subseteq \eta,  
  $$
  so $ss^*\eta\subseteq\eta$, by \lcite{\CoolLemma} and hence
$\reg_{s^*}(\eta)$ is a filter by \lcite{\MovingFilters.i}.

We claim that $\xi\subseteq \reg_{s^*}(\eta)$.  In order to prove it
let $t\in\xi$.  Then $st\in s\xi\subseteq\eta$, and hence by
\lcite{\MovingFilters.iii} we deduce that $s^*st\in \reg_{s^*}(\eta)$.
Since $t\geq s^*st$, we conclude that $t\in \reg_{s^*}(\eta)$.  This
proves our claim and since $\xi$ is an ultrafilter, we actually get
$\xi= \reg_{s^*}(\eta)$.  Therefore
  $$
  \reg_s(\xi)= \reg_s\big(\reg_{s^*}(\eta)\big) = \eta.
  \proofend
  $$

By \lcite{\MovingUltraFilters} we have that $\reg_s$ defines a map 
  $$
  \reg_s: \U_{s^*s} \to \U_{ss^*}
  $$
  which is bijective by \lcite{\VaiEVolta}.  Obviously
$\reg_s\inv=\reg_{s^*}$.

It is our next short term
goal to show that $\reg$ is a semigroup homomorphism of $\S$ into
$\I(\U)$.   The next result will be useful to help us understand the
domain of the composition of these maps.

\state Lemma
  For every $s\in\S$ and $e\in E(\S)$ one has that
$\reg_s(\U_{es^*s}) = \U_{ses^*}$.

\proof
  Let $\xi\in \U_{es^*s}=\U_e\cap\U_{s^*s}$.  Given $t\in\xi$ we have
that $et\in\xi$, and hence $set\in s\xi\subseteq\reg_s(\xi)$.  Since 
  $$
  set = ss^*set = ses^*st,
  $$
  we deduce from \lcite{\CoolLemma} that
$ses^*\reg_s(\xi)\subseteq\reg_s(\xi)$, and hence
$\reg_s(\xi)\in\U_{ses^*}$.  This shows that $\reg_s(\U_{es^*s})
\subseteq \U_{ses^*}$.  Since $ses^* = ses^*\ ss^*$, we may apply the
part of the result already proved, with $s^*$ replacing $s$, and
$ses^*$ replacing $e$, to obtain
  $$
  \reg_{s^*}(\U_{ses^*}) \subseteq \U_{s^*ses^*s} = \U_{es^*s},
  $$
  therefore
  $$
  \U_{ses^*} = 
  \reg_s\big(\reg_{s^*}(\U_{ses^*})\big) \subseteq
  \reg_s(\U_{es^*s}),
  $$
  concluding the proof.
  \proofend

From now on we will regard the $\reg_s$ as  partially defined
bijections on $\U$.  If  $f$ and $g$ are partial bijections on a set
$X$, say 
  $$
  f:A\to B \and g: C\to D,
  $$
  where $A,B,C$ and $D$ are subsets of $X$, then the
composition $gf$ is defined on $f\inv(C\cap B)$ by the expression
$fg(x) = f(g(x))$.

\state Proposition
  \label ComposingSigmaS
  For every $t,s\in\S$ one has that $\reg_t\reg_s = \reg_{ts}$.

\proof
  By the above remark the domain of the composition $\reg_t\reg_s$ is 
  $$
  \reg_s\inv(\U_{t^*t}\cap\U_{ss^*})= 
  \reg_{s^*}(\U_{t^*tss^*})= 
  \U_{s^*t^*ts} =
  \U_{(ts)^*ts},
  $$
  which coincides with the domain of $\reg_{ts}$.  Moreover for
every $\xi\in \U_{(ts)^*ts}$ we have by \lcite{\MovingFilters.iii}
that
  $$
  ts\xi =
  t(s\xi) \subseteq
  t\reg_s(\xi) \subseteq
  \reg_t(\reg_s(\xi)),
  $$
  so $\reg_t(\reg_s(\xi)) = \reg_{ts}(\xi)$, by
\lcite{\UsefulCharac}.
  \proofend
  
The following is one of our main results:

\state  Theorem 
  \label TheRegRepp
  Let $\S$ be an inverse semigroup with zero.  Then the correspondence
$s\mapsto\reg_s$ is a tight representation of $\S$ in the Boolean
inverse semigroup $\I(\U)$.

\proof
  That $\reg$ is a semigroup homomorphism follows from
\lcite{\ComposingSigmaS}, so it suffices to prove that the restriction
of $\reg$ to $E(\S)$ is a tight representation of the latter in
$E(\I(\B)) = \P(\B)$.

  If $s=0$, then $\U_{s^*s}=\U_{ss^*}=\emptyset$, and hence $\reg_s$
is the empty function, namely the zero element of $\I(\U)$, proving
\lcite{\DefineRepSemilat.i}.  As for \lcite{\DefineRepSemilat.ii} it
immediately follows from the fact that $\reg$ is multiplicative.

In order to prove tightness we would like to use
\lcite{\AltLatTightRep}, so we first need to check the validity of
either (i) or (ii) in \lcite{\AltLatTightRep}.  Thus, suppose that
\lcite{\AltLatTightRep.ii} fails, meaning that $\E(\S)$ admits a
finite cover,  say $Z$.  We will prove that 
  $$
  \U = \union\nolimits_{z\in Z}\U_z.
  $$
  By way of contradiction assume that
$\xi$ is an ultrafilter which is not in any $\U_z$.  By
\lcite{\UltraCool} for each $z\in Z$ there exists some $t_z\in\xi$
such that $zt_z=0$.  Using \lcite{\DefineFilter.iii} pick some
$t\in\xi$ such that $t_z\geq t$,  for all $z\in Z$, and notice that
  $$
  zt \leq zt_z=0, 
  $$
  so $zt=0$, and hence $ztt^*=0$, 
  which means that $z\perp tt^*$.  But since $Z$ is a cover for
$\E(\S)$ this implies that $tt^*=0$, and hence that $t=0$,
contradicting the fact that $t\in\xi$.  This proves 
\lcite{\AltLatTightRep.i}, so we may use the simplified test given
there to prove that $\lambda$ is tight.  

We therefore let $x\in \E(\S)$ be a nonzero element and 
$Z$ be a cover for $[0, x]$.
We must prove
that
  $$
  \union_{z\in Z} \U_z \supseteq 
  \U_x.
  $$
  So let $\xi$ be an ultrafilter in $\U_x$ and suppose by
contradiction that $\xi\notin\U_z$, for any $z\in Z$.  Then, by
\lcite{\UltraCool}, for each $z$ in $Z$ there exists some $t_z\in\xi$
such that
  $
  zt_z=0, 
  $
  and  by \lcite{\DefineFilter.iii} we may pick $t\in\xi$ such that
$t_z\geq t$, for all $z\in Z$.  As above this gives 
  $
  zt \leq zt_z = 0, 
  $
  so 
  $$
  zt=0 \for z \in Z.
  $$
  Given that $\xi\in\U_x$, we have that $x\xi\c\xi$,  so in particular   
$xt\in\xi$.  Now let $e = tt^*x$, and observe  that $e\leq x$, so
$e\in[0, x]$.
Moreover $e\neq0$, because 
  $$
  et =
  tt^*xt = 
  xtt^*t = 
  xt \in  
  \xi.
  $$
  By hypothesis we deduce that there exists some $z\in Z$ such that
$z\its e$,  whence
  $$
  0 \neq ze = ztt^*x = 0, 
  $$
  a contradiction.  This shows that $\xi\in\U_z$, for some $z$ in $Z$,
hence concluding the proof that $\reg$ is tight.
  \proofend

\definition We shall say that the above representation $\reg$ is the
\stress{regular representation} of $\S$.

\section{Faithfulness of tight representation}
  As in the previous section we fix an inverse semigroup $\S$ with
zero.  In this section we would like to study conditions under which
the regular representation of $\S$ is injective.  As we shall see,
injectivity does not always hold and in fact it is often the case that
different elements $s$ and $t$ in $\S$ are not separated by any tight
representation of $\S$ whatsoever.  Among our goals in this section we
will characterize precisely when does this happen.

To ease our task we will make an important assumption about the
inverse semigroup involved, which fortunately does not rule out some
important classes of inverse semigroups, such as the $E^*$-unitary
ones.

\definition
  \label DefineContinuous
  \izitem \zitem
  Let $s,t\in S$. We shall that say $s$ \stress{essentially coincides
with} $t$, in symbols $s\equiv t$, if $s^*s=t^*t$, and for every
nonzero idempotent $f\leq s^*s$, there exists a nonzero idempotent
$e\leq f$, such that $se=te$.
  \zitem We shall say that $\S$ is \stress{continuous} if $s\equiv t$
implies that $s=t$.

The fact that $s\equiv t$ is to be interpreted somewhat in the same
way as when two functions agree on a dense subset of their common
domain.  So much so that we have:

\state Proposition
  \label Localizations
  Let $\S$ be a \stress{localization} in the sense of Kumjian 
\cite{\Kumjian}, that is, $\S$ is an inverse subsemigroup of\/
$\I(X)$, where $X$ is a topological space, and $\S$ consists of
homeomorphisms between open subsets of $X$, the domains of which form
a basis for the topology of $X$.  We suppose in addition that $\S$
contains the empty function $\emptyset$, and hence $\S$ is an inverse
semigroup with zero.  Then $\S$ is continuous in the sense of
\lcite{\DefineContinuous}.

\proof
  Let $s,t\in\S$ be such that $s\equiv t$.  Identifying idempotents
with their domains, as usual, let $U=s^*s=t^*t$, and put 
  $$
  D=\{x\in U: s(x)=t(x)\}.
  $$
  We claim that $D$ is dense in $U$.  To prove it let
$A\subseteq U$ be a nonempty open set.  By hypothesis there exists a
nonzero idempotent (i.e.~an open set) $f$ such that $f\subseteq A$,
and consequently $f\leq s^*s$.  Since $s\equiv t$, we may find a
nonzero idempotent $e\leq f$, such that $se=te$.  Picking any $x\in
e$, we then have that $x\in f\subseteq A$, and $s(x) = se(x) =
te(x)=t(x)$, so $x\in A\cap D$, proving that $D$ is dense.  Since $s$
and $t$ are continuous we deduce that $s=t$.
  \proofend

Not all inverse semigroups are continuous.  Suppose for example that
$\S$ is an inverse semigroup with zero and consider
$\S'=\S\mathrel{\dot \cup}\{z\}$, where $z\notin\S$. Define a
multiplication operation on $\S'$ extending that of $\S$ and such that 
  $$
  zs=sz=z
  \for s\in\S'.
  $$
  It turns out that $\S'$ is an inverse semigroup with zero, except
that the \stress{zero} of $\S'$ is $z$, rather than the original zero
of $\S$ (which we denote by $0$).

  Given any $s,t\in\S$ with $s^*s=t^*t$, and any nonzero (i.e.~different
from $z$) idempotent $f\leq s^*s$, notice that 0 is a nonzero (sic)
idempotent with $0\leq f$, and $s0=0=t0$.  Thus $s\equiv t$ even though
$s$ and $t$ might not coincide.

The following additional counter-example is due to Szendrei (personal
communication).  Let $X=\{1, 2, 3\}$, and let $\S$ be the inverse
subsemigroup of $\I(X)$ consisting of the following four elements:

  \bitem $1$ -- identity permutation 
  \bitem $0$ -- empty mapping (the zero element)
  \bitem $i$ -- the partial identity sending 1 to 1 and undefined
otherwise
  \bitem $s$ -- the transposition interchanging 2 and 3 (and sending 1
to 1).

\medskip\noindent Clearly $1$, $0$, and $i$ are idempotents of $\S$
forming a three-element chain.  We have $s^*s = 1$, and for both
nonzero idempotents $e$ with $e \leq s^*s$ (that is, for both $e =
1$ and $e = i$), the relations $i \leq e$ and $si = 1i$ hold,  whence
$s\equiv 1$.

Recall that an inverse semigroup with zero is said to be $E^*$-unitary
\cite{\Szendrei}, \scite{\Lawson}{Section 9}, if whenever $s\in\S$,
and $se=e$, for some nonzero idempotent $e$, then $s$ is necessarily
also idempotent.

\state Proposition 
  \label EsUnitaryCont
  Every $E^*$-unitary inverse semigroup with zero is
continuous.

\proof 
  Let $\S$ be an $E^*$-unitary inverse semigroup with zero and let
$s,t\in\S$ be such that $s\equiv t$.  Plugging $f=s^*s$ in the
definition there exists a nonzero idempotent $e\leq s^*s$, such that
$se=te$.  Then \scite{\actions}{5.3} applies giving $s=t$.
  \proofend

We now return to studying the general case.

\state Proposition 
  \label UnfaithfulButRelated
  Let $\S$ be an inverse semigroup with zero.
  If $s,t\in\S$ are such that $s^*s=t^*t$, and $\reg_s=\reg_t$, where
$\reg$ is the regular representation of $\S$, then $s\equiv t$.

  \proof
  Given a nonzero idempotent $f\leq s^*s$, choose an ultrafilter $\xi$
such that $f\in\xi$.  Since $s^*sf = f \in\xi$, we have by
\lcite{\CoolLemma} that $\xi\in\U_{s^*s}$.  Thus $\reg_s(\xi) =
\reg_t(\xi)$.  In addition we have by \lcite{\MovingFilters.iii} that
  $$
  sf\in \reg_s(\xi) = \reg_t(\xi)\ni tf,
  $$
  so by \lcite{\DefineFilter.iii} there exists $u\in\reg_s(\xi)$
such that $sf,tf\geq u$.  Therefore 
  $$
  sfu^*u = u = tfu^*u.
  $$
  Thus $e:= fu^*u$ \ is a nonzero idempotent (because $u\neq0$) such
that $e\leq f$, and $se=te$.  This proves that $s\equiv t$.
\proofend

Although it is not crucial for our purposes it would be interesting to
decide if the converse of the above result holds.

In the following  main result we characterize precisely the extent to
which tight representations do not separate points of a continuous
inverse semigroup $\S$.

\state Theorem
  \label HowNotFaithful
  Let $\S$ be a continuous inverse semigroup with zero and let
$s,t\in\S$.  Then the following are equivalent:
  \izitem
  \zitem $\sigma(s) = \sigma(t)$ for every tight representation
$\sigma$ of\/ $\S$,
  \zitem $\reg(s)=\reg(t)$,
  \zitem $st^*t = ts^*s$, and $s^*st^*t$ is dense in both $s^*s$ and $t^*t$,
  \zitem $tt^*s = ss^*t$, and $ss^*tt^*$ is dense in both $ss^*$ and $tt^*$.

\proof
  (i) $\imply$ (ii): \  obvious.

\medskip\noindent (ii) $\imply$ (iii): \ 
  If $\reg(s)=\reg(t)$ then in particular the domains of $\reg(s)$ and
$\reg(t)$ must coincide, and hence $\U_{s^*s} = \U_{t^*t}$.  The last
assertion in (iii) then follows from \lcite{\SameUs}.
  Next let $\alt s= st^*t$, and $\alt t=ts^*s$, and observe that 
  $$
  \alt t^*\alt t = s^* st^*t = \alt s^* \alt s.
  $$
  Moreover we have
  $$
  \reg(\alt s) = \reg(st^*t) =
  \reg(s) \reg(t)\inv \reg(t) = 
  \reg(t) \reg(s)\inv \reg(s) = 
  \reg(ts^*s) = 
  \reg(\alt t).
  $$
  Invoking \lcite{\UnfaithfulButRelated} we conclude that $\alt
s\equiv \alt t$, and hence that $\alt s= \alt t$, because $\S$ is
continuous.

\medskip\noindent (iii) $\imply$ (i): \ 
  Let $\sigma$ be a tight representation of $\S$.  
  Since $s^*st^*t$ is dense in both $s^*s$ and $t^*t$, we have by the
last part of \lcite{\SameUs} that 
$\sigma(s^*s) = \sigma(t^*t)$.  Therefore
  $$
  \sigma(s) =
  \sigma(s s^*s) =
  \sigma(s) \sigma(s^*s) =
  \sigma(s) \sigma(t^*t) = 
  \sigma(st^*t) \$= 
  \sigma(ts^*s) =
  \sigma(t)\sigma(s^*s) = 
  \sigma(t)\sigma(t^*t) = 
  \sigma(tt^*t) = 
  \sigma(t).
  $$

\medskip\noindent (ii) $\Leftrightarrow$ (iv): \ 
  Since $\reg(s^*) = \reg(s)\inv$, and similarly for $t$, we have that
(ii) is equivalent to saying that $\reg(s^*)=\reg(t^*)$.  Exchanging
$s$ and $t$, respectively by $s^*$ and $t^*$, and applying the already
proved equivalence between (ii) and (iii), we then see that (ii) is
equivalent to saying that $s^*tt^* = t^*ss^*$, and that $ss^*tt^*$ is
dense in both $ss^*$ and $tt^*$, which is tantamount to (iv).
  \proofend

\section{A counter-example}
  Given the use of the continuity hypothesis in the proof of the
implication (ii) $\imply$ (iii) of \lcite{\HowNotFaithful} it is
interesting to decide whether or not that result survives
in the absence of such a hypothesis.  In this section we present an
example to show that it does not.

We begin by exhibiting a general construction of inverse
semigroups. In order to do so recall from \qcite{\PatBook}{?} that a
\stress{congruence} on an inverse semigroup $\S$ is an equivalence
relation ``$\sim$" such that $us\sim ut$ and $su\sim tu$, whenever
$s\sim t$ and $u\in\S$.  Given any such relation the quotient set \
$\S/{\sim}$ \ is an inverse semigroup \qcite{\PatBook}{?}.  Our construction
will be attained by means of taking a quotient.

Let $\E$ be a semilattice with smallest element 0, and let $G$ be a
group.  Suppose that for each $x\in\E$ we are given a normal subgroup
  $$
  N_x\trianglelefteq G, 
  $$
  such that $N_0=G$, and whenever $x\leq y$ in $\E$ one has that
$N_x\supseteq N_y$.

Viewing both $\E$ and $G$ as inverse semigroups, consider their
cartesian product
  $$
  \T = \E\times G,
  $$
  with coordinatewise operations.  Clearly $\T$ is an inverse
semigroups as well.  We define a congruence on $\T$ by saying that
  $$
  (x,g) \sim (y,h),
  $$
  if and only if $x=y$ and $h\inv g\in N_x$.  

If $\pi_x$ denotes the
quotient map from $G$ to $G/N_x$, then the last condition above is
perhaps more conveniently stated by saying that $\pi_x(g)=\pi_x(h)$.
By our assumptions about the $N_x$ it is evident that, whenever $x\leq
y$, one has that
  $$
  \pi_y(g)=\pi_y(h) \ \imply \  \pi_x(g)=\pi_x(h) 
  \for g,h\in G.
  \eqmark BiggerRel
  $$

  Leaving aside the obvious verification that ``$\sim$" is an
equivalence relation, let us check that it is indeed a congruence.
For this suppose that 
  $(x,g) \sim (y,h)$, 
  and let $(z,k)\in\T$.  Then 
  $$
  (x,g)(z,k) = (x\inf z, gk) \and
  (y,h)(z,k) = (y\inf z, hk),
  $$
  and we must prove that $(x\inf z, gk) \sim (y\inf z, hk)$.
Obviously $x=y$, so $x\inf z = y\inf z$, as well.  It then suffices to
prove that
  $\pi_{x\inf z}(gk) = \pi_{x\inf z}(hk)$.  Noticing that
$\pi_x(g)=\pi_x(h)$, we have
  $$
  \pi_{x\inf z}(gk) = 
  \pi_{x\inf z}(g) \ \pi_{x\inf z}(k) \={(\BiggerRel)}
  \pi_{x\inf z}(h) \ \pi_{x\inf z}(k) =
  \pi_{x\inf z}(hk).
  $$
  The proof that ``$\sim$" is invariant under right multiplication is
done in a similar way.

  \def\eq#1#2{[#1, #2]}

  We shall let $\S=\S(\E,G,\{N_x\}_x)$ be the inverse semigroup
obtained by taking the quotient of $\T$ by ``$\sim$".  Given
$(x,g)\in\T$, we will henceforth refer to its equivalence class by
$\eq xg$.  Notice that for every $(x,g)\in\T$ one has that 
  $$
  \eq 01 \eq xg = \eq{0\inf x}{g} = \eq{0}{g} = \eq{0}{1},
  $$
  the last equality following from the assumption that $N_0=G$.
One may similarly prove that $\eq xg   \eq 01 = \eq{0}{1}$, which
means that $\eq 01$ is a zero element for $\S$.  
Moreover
  $$
  \eq xg 
  {\eq xg}^*
  = \eq xg \eq x{g\inv} = \eq x1,
  $$
  so $E(\S)$ consists of the set of all equivalence classes $\eq x1$,
for $x\in \E$.  Since
  $
  (x,1)\sim (y,1)
  $
  if and only if $x=y$, we deduce that $E(\S)$ is isomorphic to $\E$.

We will now consider a more concrete application of these ideas.  Let
$\{0,1\}$ have the obvious Boolean algebra structure and put
$\B=\{0,1\}\times\{0,1\}$.  The unit of $\B$ is clearly $1 = (1,1)$,
and its zero element is $0 = (0,0)$.  We shall denote the remaining
elements by
  $$
  e_1 = (1,0) \and e_2 = (0,1),
  $$
  so $\B=\{0,e_1,e_2,1\}$.
We will temporarily view $\B$ simply as a semilattice.
Given the sheer simplicity of $\B$, given any two
normal subgroups $N_1,N_2\trianglelefteq G$, and setting 
  $$
  \eqmatrix{
  N_0 & = & G, \hfill \cr
  N_{e_1} & = & N_1, \hfill \cr
  N_{e_2} & = & N_2,  \hfill \cr
  N_1 & = & \{1\}, \hfill  \cr
  }
  $$
  one may easily check that the collection $\{N_x\}_{x\in\B}$
satisfies the conditions above so that  we may construct the
associated inverse semigroup $\S=\S(\B,G,\{N_x\}_x)$ as above.  

 As already noticed $E(\S) = \B$, so $\S$ is a Boolean inverse
semigroup.  By \lcite{\BooleToBoole} one sees that the identity mapping
  \ $
  \iota: \S \to \S
  $ \
  is a tight representation of $\S$ in itself.  It follows that
\lcite{\HowNotFaithful.i} only holds when $s=t$.  However, we will
show that \lcite{\HowNotFaithful.ii} might hold for $s\neq t$.

Before we begin let us agree on a particularly useful notation for
elements of $\S$.  Noticing that $N_1=\{1\}$ observe that 
$(1,g)\sim(1,h)$ if and only if $g=h$.  Thus the mapping 
  $$
  g\in G\mapsto \eq 1g\in\S
  $$
  is a semigroup monomorphism, and hence we may identify $G$ with its
copy within $\S$.  We have already observed that 
  $$
  x\in \B\mapsto \eq x1\in\S
  $$
  is also a semigroup monomorphism and hence we are allowed to think
of $\B$ as a subsemigroup of $\S$. 
  Given any $(x,g)\in\T$ we have that
  $$
  \eq xg = \eq x1 \eq 1g = xg,
  $$
  where in the last term we are fully enforcing our identifications.
Therefore $\S = \B G$, and thanks to the product structure of $\T$ notice
moreover that  $\B$ and $G$ commute.  

With the purpose of understanding the order structure of $\S$ let $\eq
xh$ and $\eq yg$ be elements in $\S$ with $\eq xh\leq\eq yg$,
that is,
  $$
  \eq xh =
  \eq yg \eq xh \eq xh^* =
  \eq yg \eq x1 =
  \eq {y\wedge x}g.
  $$
  This is the same as saying that $x\leq y$, and $\pi_{x}(h)
=\pi_x(g)$.  This implies in particular that $\eq xh=\eq xg$,  so
every inequality in $\S$ is of the form
  $$
  \eq xg\leq\eq yg,
  $$
  with $x\leq y$.  The only nontrivial inequalities (i.e.~not
involving zero nor an equality) are therefore
  $$
  ge_i\leq g,
  $$
  for $i=1,2$, and $g\in G$.  It follows that the minimal nonzero
elements of $\S$ are precisely those of the form $ge_i$, as above.
This said it is easy to see that the most general ultrafilter in $\S$
is
  $$
  \xi_{ge_i} = \{s\in\S: s\geq ge_i\} = \{ge_i, g\},
  $$
  for $i=1,2$, and $g\in G$.  

For the purpose of giving our counter-example we will suppose in
addition that $N_{e_1}\cap N_{e_2}\neq \{1\}$.  Given $s\in G$, choose
a nontrivial element $n\in N_{e_1}\cap N_{e_2}$ and put $t=sn$.
Clearly $s\neq t$, but
  $$
  \pi_{e_i}(s)=\pi_{e_i}(t)
  \for i=1,2.
  $$
  It is our intention to prove that $\reg_{s}=\reg_{t}$.  In order to
do so notice that
$s^*s=t^*t=1$, and that $\U_1=\U$,  so the domain of both $\reg_s$ and
$\reg_t$  is the set of all ultrafilters.  

Given any ultrafilter $\xi$, write 
$\xi=\xi_{ge_i}$,  for some $g\in G$, and $i=1, 2$.
  Since $ge_i\in\xi$, we have that $sge_i\in\reg_s(\xi)$, by
\lcite{\MovingFilters.iii}.  Recalling that $sge_i$ is a minimal
element, we necessarily have
  $
  \reg_s(\xi) = \xi_{sge_i},
  $
  and similarly $\reg_t(\xi) = \xi_{tge_i}$.  Moreover
  $$
  (e_i,sg)\sim (e_i,tg),
  $$
  because $\pi_{e_i}(sg) = \pi_{e_i}(tg)$, so $sge_i=tge_i$, and hence
  $$
  \reg_s(\xi) = \xi_{sge_i} = \xi_{tge_i} = \reg_t(\xi).
  $$
  Since $\xi$ is arbitrary we deduce that $\reg_s= \reg_t$.  The big
conclusion is that \lcite{\HowNotFaithful.ii} holds for $s$ and $t$,
but \lcite{\HowNotFaithful.i} does not.  The trouble is of course that
$\S$ is not continuous, and this concludes our goal of showing that
\lcite{\HowNotFaithful} cannot be proved without the continuity
hypothesis.

  \endgroup

  \begingroup
  \bigskip\bigskip 
  \font \sc = cmcsc8 \sc
  \parskip = -1pt

  Departamento de Matem\'atica 

  Universidade Federal de Santa Catarina

  88040-900 -- Florian\'opolis -- Brasil

  \eightrm r@exel.com.br

  \endgroup
  \end